\documentclass[12pt]{article}
\usepackage{amssymb,amsmath,amstext, amsthm}

\newtheorem*{main}{Main Theorem}
\newtheorem{theorem}{Theorem}[section]
\newtheorem{corollary}[theorem]{Corollary}
\newtheorem{proposition}[theorem]{Proposition}
\newtheorem{lemma}[theorem]{Lemma}

\theoremstyle{definition}
\newtheorem{definition}[theorem]{Definition}

\renewcommand{\sharp}{\#}
\def\al{\alpha}
\def\be{\beta}
\def\ga{\gamma}
\def\de{\delta}
\def\ep{\epsilon}
\def\la{\lambda}
\def\Om{\Omega}
\def\La{\Lambda}

\def\P{{\mathbb P}}
\def\R{{\mathbb R}}

\def\C{{\mathbb C}}
\def\S{S}
\def\G{G}

\def\CF{{\cal F}}
\def\CC{{\cal C}}
\def\CH{{\cal H}}
\def\CB{{\cal B}}

\def\CQ{{\cal Q}}
\def\CT{{\cal T}}

\def\bc{\mathbf{c}}
\def\bcdot{\mathbf{\dot c}}
\def\bcdotdot{\mathbf{\ddot c}}

\def\e{\mathrm{{e}}}

\newcommand{\ra}{\rightarrow}

\newcommand{\binomial}[2]{\left(\begin{matrix}#1\\#2\end{matrix}\right)}
\newcommand{\primal}{\ensuremath{\mathbf{P}}}
\newcommand{\dual}{\ensuremath{\mathbf{D}}}
\newcommand{\primaldual}{\ensuremath{\mathbf{PD}}}
\newcommand{\im}{\mbox{im }}
\newcommand{\da}{\dagger}
\newcommand{\rank}{\mbox{Rank }}
\newcommand{\diag}{\mbox{Diag }}
\newcommand{\codim}{\mbox{codim }}

\title{On the Curvature of the Central Path of Linear Programming Theory.%
\thanks{Mathematics Subject Classification (MSC2000): 90C51 (Primary),
90C60, 68Q25, 65H10 (Secondary).}}
\author{Jean-Pierre Dedieu%
        \thanks{%
               MIP. D\'epartement de Math\'ematique,
               Universit\'e Paul Sabatier,
               31062 Toulouse cedex 04, France
                ({\tt dedieu@mip.ups-tlse.fr}).
}
\and    Gregorio Malajovich%
        \thanks{
                Departamento de Matem\'atica Aplicada,
                Universidade Federal de Rio de Janeiro,
                Caixa Postal 68530,
                CEP 21945-970, Rio de Janeiro, RJ, Brazil
                ({\tt gregorio@ufrj.br}) Gregorio Malajovich
                was partially sponsored by CNPq and by FAPERJ
        (Brazil). His visits
                to Universit\'e Paul Sabatier were supported by the
                International Cooperation Agreement Brazil-France
                in Mathematics.
}
\and    Mike Shub%
        \thanks{%
                Department of Mathematical Sciences,
                T. J. Watson Research Center,
                Yorktown Heights, NY 10598, USA,
                and:
                Department of Mathematics,
                University of Toronto,
                100 St. George Street,
                Toronto, Ontario M5S 3G3, Canada
                ({\tt shub@math.toronto.edu}).
}}

\begin{document}
\maketitle

\vskip 1.5cm \begin{center}{\bf We dedicate this paper with great
admiration and affection\\ to our friend and teacher Steve Smale
for his seventy fifth birthday.}\end{center} \vskip 1.5cm
\begin{abstract}

We prove a linear bound on the average total curvature of the
central path of linear programming theory in terms on the number
of independent variables of the primal problem, and independent of
the number of constraints.
\end{abstract}

\section{Introduction.}
\label{intro}

Consider a linear programming problem in the following
primal/dual form:

\[
\min_{
\begin{array}{c}
Ax - s = b \cr s \ge 0 \cr
\end{array}}
\left\langle c,x \right\rangle
\hspace{2em} \text{and} \hspace{2em}
\max_{
\begin{array}{c}
A^Ty = c \cr y \ge 0 \cr
\end{array}}
\left\langle b,y \right\rangle
\ .
\]

Here $m>n \ge 1$ and $A$ is an $m \times n$ real matrix assumed to
have rank $n$,  $b \in \R^m$ and $c \in \R^n$ are given vectors
and $c$ is non-zero, $y, s \in \R^m$ and $x \in \R^n$ are unknown
vectors ($s$ is the vector of slack variables).

Our principal result bounds the total curvature of the union of
all the central paths associated {with} all the feasible
regions obtained by considering all the $2^m$ possible sign
conditions
\[
s_i \ \ep_i \ 0, \hspace{2em} i=1,\cdots,m,
\]
where $\ep_i$ is either $\geq$ or $\leq.$

Formal definitions will be given in subsequent sections. The rest
of the results in the introduction follow from the next theorem
which requires the rest of the paper.

\begin{theorem}\label{th-main-sum} 
Let $m > n \ge 1$. Let $A$ be an $m \times n$ matrix of rank $n$,
and let $b \in \R^m$ and $c \in \R^n$, $c$ non-zero. The sum over
all $2^m$ sign conditions of the total curvature of the
primal/dual central paths (resp. primal central paths, dual
central paths) is less than or equal to $2 \pi n
\binomial{m-1}{n}$ (resp. $2 \pi (n-1) \binomial{m-1}{n}$, $2 \pi
n \binomial{m-1}{n} .$)
\end{theorem}

Theorem \ref{th-main-sum} allows us to conclude various results on
the average curvature of the central paths corresponding to
various probability measures on the space of problems. We begin
with our main motivating example.

Central paths are numerically followed to the optimal solution of
linear programming problems by interior point methods. For
relevant background material on interior point methods see Renegar
\cite{ren}. Our point in studying the total curvature is that
curves with small total curvature may be easy to approximate with
straight lines. So, small total curvature may contribute to the
understanding of why long step interior point methods are seen to
be efficient in practice. In Dedieu-Shub \cite{ded} we studied the
central paths of linear programming problems defined on strictly
feasible compact polyhedra (polytopes)\footnote[1]{A feasible
region for a linear programming problem is a polyhedron, a compact
polyhedron is a polytope. 
}  
from a dynamical systems perspective. In this
paper we optimistically conjectured that the worst case total
curvature of a central path is $O(n)$. Our first average result
and main theorem lends some credence to this conjecture, proving
it on the average.

If we assume that the {primal polyhedron $\{x|Ax - b = s
\ge 0\}$} is compact and strictly feasible {(i.e., has
non-empty interior)}, then the primal and dual problems have
central paths which are each the projection of a primal/dual
central path and all these central paths lead to optimal
solutions. So for our purposes we will get a meaningful number if
we divide the total curvature of the central paths of the all the
strictly feasible polytopes arising from all possible
sign conditions by the number of distinct strictly
feasible polytopes associated with the $2^m$ sign conditions:
\[
Ax - s = b \hspace{2em} s_i \ \ep_i \ 0 \hspace{1em} i=1, \dots,
m,
\]
where $\ep_i$ is either $\ge$ or $\le$. The cardinality of the set
of these polytopes is $\leq  \binomial{m-1}{n}$ and equality holds
for almost all $(A,b)$, see section \ref{formal}. When equality
holds we say $(A,b)$ is in $\textbf{general position}.$.

We use Theorem \ref{th-main-sum} to give an upper bound on the sum
of the curvatures.

We obtain the following average result.

\begin{main}\label{th-main-informal}
Let $m > n \ge 1$. Let $A$ be an $m \times n$ matrix of rank $n$,
and let $b \in \R^m$ and $c \in \R^n$, $c$ non-zero such that
(A,b) is in general position. Then the average total curvature of
the primal/dual central paths (resp. primal central paths, dual
central paths) of the strictly feasible polytopes defined
by (A,b) is less than or equal to $2 \pi n$ {(resp. $2
\pi (n-1)$, $2 \pi n .$)}
\end{main}


We may also average over more general probability measures on the
data $A$, $b$, $c$ defining the problem. First we more precisely
define the space of problems $\emph{P}$ and measures $\mu,\nu$ we
consider. $\emph{P}=\emph{I}\times
\mathbb{R}^m\times\mathbb{R}^n.$ Here $\emph{I}$ is the open set
of $rank(n)$, $m$ by $n$ real matrices,and we assume for
convenience that no row of any element of $\emph{I}$ is
identically zero. Let $\emph{D}$ be the group with $2^m$ elements
consisting of those $m$ by $m$ diagonal matrices whose diagonal
entries are all either $1$ or $-1$. So for $D \in \emph{D}$, D
acts on $\emph{P}$ by $D((A,b,c))=(DA,Db,c).$ The set of problems
defined by the orbit of $(A,b,c)$ under the action of $\emph{D}$
is the same as considering $(A,b,c)$ with all possible sign
conditions, so each orbit has $2^m$ distinct elements. We say that
a probability measure $\mu$ is \textbf{sign invariant} if it is
invariant under the action of $\emph{D}$, i.e. $D_*\mu = \mu$ for
all $D \in \emph{D}$.

We now generalize Theorem \ref {th-main-informal} once again
averaging over problems with a strictly feasible primal
polytope.

Let $\mu$ be a sign-invariant probability measure on the data $A$,
$b$, $c$. If the the set of $(A,b,c)$ in $\emph{P}$ such that
$(A,b)$ are in general position has full measure we will say that
$\mu$ is \textbf{full (for general position)}. This is the case
for example if $\mu$ is supported on a finite union of orbits of
$\emph{D}$ through elements in general position or if $\mu$ is
absolutely continuous with respect to the Lebesgue measure, see
section \ref{formal}. For instance, an independent Gaussian
probability distribution with zero average and arbitrary variance
for each coefficient of the data is sign invariant and full.

\begin{corollary}\label{cor-main-av}
  Let $m>n$ and let $\mu$ be a sign-invariant and full (for general
  position) probability measure on $\emph{P}$. Let $\mathrm{Feas}$ be the set
of data $A,b,c$ with a strictly feasible primal polytope.
Let $\nu$ be the conditional probability measure (with respect to
$\mathrm{Feas}$) defined for any measurable $V$ by
\[
\nu (V) = \frac{ \mu(V \cap \mathrm{Feas}) }{ \mu(\mathrm{Feas})} .
\]
  Then, the average (w.r.t. $\nu$) total curvature of the
primal/dual central path (resp. primal central path, dual central path)
is less than or equal to $2 \pi n$ (resp. $2 \pi (n-1)$, $2 \pi n .$)
\end{corollary}

This corollary while almost immediate requires a little proof
which we carry out in section \ref{formal}. There is another
version of Corollary \ref{cor-main-av} which is perhaps a little
more natural form the point of view of regions which have central
paths defined for all positive parameter values. We state it below
but don't prove it as the proof is the same as for Corollary
\ref{cor-main-av}. For a primal/dual central path to exist for all
positive parameter values a necessary and sufficient condition is
that both primal and dual problems are strictly feasible see
\cite{tod},\cite{wri} . If this is the case we say that the
primal/dual polyhedra are \textbf{jointly strictly feasible}.
Every strictly feasible primal polytope gives rise to
primal/dual jointly strictly feasible polyhedra, but there are
more of the later generally among the polyhedra arising from the
$2^m$ possible sign conditions in a linear programming problem.
Generally the number of jointly strictly feasible primal/dual
polyhedra is $\binomial{m}{n}$. We may see this simply since there
are generally  $\binomial{m}{n}$ vertices to the primal polyhedra
and at each vertex almost all non-zero $c$ select a unique primal
polyhedron for which that vertex minimizes the optimization
problem, see\cite{adl}. When the number is $\binomial{m}{n}$ we
say that $(A,b,c)$ is in $\textbf{joint general position}$. If we
consider a sign invariant probability measure which is full (for
joint general position) i.e. the set of problems $(A,b,c)$ which
are in joint general position has full measure,  we get a slight
improvement of Corollary \ref{cor-main-av}.

\begin{corollary}\label{cor-main-alt}
  Let $m>n$ and let $\mu$ be a sign-invariant and full (for joint general
  position) probability measure on $\emph{P}$. Let $\mathrm{Feas}$ be the set
of data $A,b,c$ with  joint strictly feasible primal/dual
polyhedra. Let $\nu$ be the conditional probability measure (with
respect to $\mathrm{Feas}$) defined for any measurable $V$ by
\[
\nu (V) = \frac{ \mu(V \cap \mathrm{Feas}) }{ \mu(\mathrm{Feas})}
.
\]
  Then, the average (w.r.t. $\nu$) total curvature of the
primal/dual central path (resp. primal central path, dual central
path) is less than or equal to $2 \pi n \frac{m-n}{m}$ (resp. $2
\pi (n-1)\frac{m-n}{m}$, $2 \pi n \frac{m-n}{m}.$)
\end{corollary}

\section{Description of the central path.}

When the optimal condition is attained, the
primal and dual problems have the same value and
the optimality conditions may be written as
\[
\left\{ \
\begin{array}{l}
Ax - s = b \cr A^Ty = c\cr sy = 0\cr y \ge 0 , \ s \ge 0\cr
\end{array}
\right.
\]
\noindent where $sy$ denotes the componentwise product of these two vectors.
The {\em primal/dual central path} of this problem is the curve
$\left(x(\mu),s(\mu),y(\mu)\right)$, $0 < \mu < \infty$, given by
\begin{equation}\label{centralpath}
\left\{ \
\begin{array}{l}
Ax - s = b \cr A^Ty = c\cr sy = \mu \e\cr y > 0 , \ s > 0\cr
\end{array}
\right.
\end{equation}
where $\e$ denotes the vector in $\R^m$ of all $1$'s.
\medskip

The {\em primal central path} is the curve
$\left(x(\mu),s(\mu)\right)$, $0 < \mu < \infty$, defined as the
curve of minimizers of the function $-\mu\sum_1^m \ln(s_i) + c^Tx$
restricted to the primal polyhedron. By the use of Lagrange
multipliers one sees that this is the curve defined by the
existence of a vector $y(\mu)$ satisfying the equations
(\ref{centralpath}). Thus the primal central path is the
projection of the primal/dual central path into the $(x,s)$
subspace.
\medskip

Similarly, the {\em dual central path} is the curve $y(\mu)$, $0 <
\mu < \infty$, defined as the curve of maximizers of the function
$\mu\sum_1^m \ln(y_i) + b^Ty$ restricted to the dual polyhedron.
By use of Lagrange multipliers one sees that this curve is defined
by the existence of vectors $x(\mu),s(\mu)$ satisfying
(\ref{centralpath}). So the dual central path is the projection of
the primal/dual central path on the $y$ subspace.

Note, as we have alluded to in the introduction, that
{when} the primal polyhedron is compact and strictly
feasible the primal central path is defined for all $0 < \mu <
\infty$ and then so are the primal/dual and dual central paths.

\section{Curvature.}\label{sec-curv}

Let $\bc : [a,b] \ra \R^n$ be a $C^2$ map with non-zero
derivative: $ \dot {\bc} (t) \not = 0$ for any $t \in
[a,b]$. We denote by $l$ the arc length:
$$l(t) = \int_a^t \| \dot {\bc}(\tau) \| \ d\tau .$$
To the curve ${\bc}$ is associated another curve on the
unit sphere, called the Gauss curve, defined by
$$t \in [a,b] \ra \ga(t)=\frac{\dot {\bc}(t)}{\| \dot {\bc}(t) \| } \in \S^{n-1}$$
which may also be parameterized by the arc length $l$ of $\bc$:
$$l \in [0,L] \ra {\dot {\bc}(l)} \in \S^{n-1}$$
with $L$ the length of the curve ${\bc}$. The curvature is
$$\kappa (l) =  \frac{d}{dl} \dot {\bc}(l) $$
see Spivak \cite[chapter 1]{spi}. In terms of the original parameter
we have
\begin{equation}
\label{curv-t}
\kappa (t) =  \frac{1}{\| \dot {\bc}(t) \|}\ \frac{d}{dt}
\left(\frac{\dot {\bc}(t)}{\| \dot {\bc}(t)
\|}\right) = \frac{ \ddot {\bc}(t) \ \| \dot {\bc}(t)
\|^2 - \dot {\bc}(t)\ \left\langle \dot {\bc}(t) ,
\ddot {\bc}(t) \right\rangle }{\| \dot {\bc}(t)
\|^4}.
\end{equation}

The total curvature $K$ is the integral of the norm of the
curvature vector:
$$K = \int_0^L \| \kappa (l) \| \ dl \ .$$
Thus, $K$ is equal to the length of the Gauss curve on the unit
sphere $\S^{n-1} \subset \R^n.$ To compute $K$ we use integral
geometry, the next section is devoted to that.

\section{An integral geometry formula.}

Let $\ga (t)$, $a \leq t \leq b$, be a $C^1$ parametric curve
contained into the unit sphere $\S^{n-1}$ with at most a countable
number of singularities (i.e. $\dot \ga (t) = 0$). The parameter interval is not necessarily
finite: $- \infty \leq a \leq b \leq \infty .$ Let us denote by
$\G_{n,n-1}$ the Grassmannian manifold of hyperplanes through the
origin contained in $\R^n$. We also denote by $d\G(\CH)$ the unique
probability measure on $\G_{n,n-1}$ invariant under the action of
the orthogonal group.

\begin{theorem} \label{th-length} The length of $\ga$ is equal to
$$L(\ga) = \int_a^b \left\| \frac{d}{dt} \ga(t) \right\|\ dt = \pi \int_{\CH \in \G_{n,n-1}} \sharp (\CH \cap \ga) \ \ d\G(\CH)$$ where
$\sharp (\CH \cap \ga)$ denotes the number of parameters $a \leq t
\leq b$ such that $\ga (t) \in \CH$: $ \sharp (\CH \cap \ga)$ is the number of intersections counted with multiplicity.
\end{theorem}

\proof If $\ga$ is an embedding then Theorem \ref{th-length}
follows from Santal\'o \cite[chapter 18, section 6]{san}, or also see
Shub and Smale \cite[section 4]{bez5}, where a similar theorem is
proved for projective spaces or Edelman and Kostlan \cite{ede}.
Now the set of $t$ such that $\frac{d}{dt} \ga(t) \not = 0$ may be
written as a countable union of intervals on each of which $\ga$
is an embedding.  \qed
\medskip

\begin{definition} \label{def-genericity} The {parametric} curve $\ga$ is transversal to $\CH \in \G_{n,n-1}$ (we also say $\CH$ is transversal to $\ga$) when $\dot \ga (t) \not \in \CH$ at the intersection points.
\end{definition}

\begin{corollary} \label{cor-length} If the number of intersections counted with multiplicity satisfies $\sharp (\CH \cap \ga) \leq \CB$ for all transversal $\CH \in \G_{n,n-1}$ then
$$L(\ga) = \int_a^b \left\| \frac{d}{dt} \ga(t) \right\|\ dt \leq \pi \CB.$$
\end{corollary}

\proof By a usual application of Sard's Theorem, see Golubitsky-Guillemin \cite{gol}, non-transversality is a zero measure event. Thus, the integral giving $L(\ga)$ only needs to be evaluated on the set $\CT$ of $\CH \in \G_{n,n-1}$ such that $\CH$ is transversal to $\ga$. Since $d\G(\CH)$ is a probability measure we get
$$L(\ga) = \pi \int_{\CH \in \CT} \sharp (\CH \cap \ga) \ \ d\G(\CH) \leq  \pi \CB \int_{\CH \in \CT} d\G(\CH)
= \pi \CB .$$ \qed

In order to bound the number of transversal intersections of the Gauss
curve with a hyperplane $\CH$, we will need the following
fact: let
\[
\begin{array}{llcl}
F :& \R \times \R^r & \rightarrow & \R^r \\
& (\mu, z) & \mapsto &F_{\mu}(z)
\end{array}
\]
be of class $\CC^2$, and assume that we are in
the conditions of the Implicit
Function Theorem, namely
$F_{\mu_0}(\bc_0) = 0$ and
$DF_{\mu_0} (\bc_0)$ (the derivative of $F$ with
respect to the $z$ variables)
has full rank.
Let $\bc(\mu): [\mu_0 - \epsilon, \mu_0 + \epsilon] \rightarrow \R^r$
be the associated implicit function, $\bc(\mu_0)=\bc_0$
and $F_{\mu} (\bc(\mu)) = 0$, and let $\bcdot(\mu)$ denote
the derivative of $\bc$ with respect to $\mu$.

Let $\CH$ denote a hyperplane, with normal vector $h$:
$$\CH = \{ z \in \R^r \ : \ \langle h , z \rangle = 0 \}.$$

\begin{lemma} \label{transverse} In the conditions above, if the Gauss curve
$\ga(\mu) = \bcdot (\mu)/\|\bcdot(\mu)\|$ intersects $\CH$
transversally for $\mu = \mu_0$, then $(\bc(\mu_0),
\bcdot (\mu_0), \mu_0)$
is a zero of the function
\begin{equation} \label{ift}
\Phi (\bc, \bcdot, \mu) = \left[
\begin{array}{c}
F_{\mu}( \bc ) \\
DF_{\mu}( \bc ) \ \bcdot  + \dot F_{\mu}(\bc)\\
\langle h, \bcdot \rangle
\end{array}
\right]
\end{equation}
Moreover, $D\Phi$ has full rank at that point.
\end{lemma}

\begin{proof}
Equation (\ref{ift}-1) and (\ref{ift}-2) are the Implicit
Function Theorem, and equation (\ref{ift}-3) is the
intersection hypothesis. We write $D\Phi (\bc, \bcdot, \mu)$ as the block matrix:
\begin{equation} \label{dPhi}
D\Phi = \left[
\begin{array}{ccc}
DF_{\mu}( \bc )
&0
&\dot F_{\mu}( \bc )
\\
D^2F_{\mu}( \bc ) \otimes \bcdot
+ D\dot F_{\mu} (\bc)
&
DF_{\mu}( \bc )
& D\dot F_{\mu}( \bc ) \ \bcdot
+
\ddot F_{\mu}( \bc )
\\
0& h^T & 0 \\
\end{array}
\right]
\end{equation}
where $D^2F_{\mu}( \bc ) \otimes \bcdot$ is the
linear map $y \mapsto D^2F_{\mu}( \bc ) (\bcdot,y)$.
By hypothesis, $DF_{\mu}( \bc )$ is invertible. Hence, the
block LU factorization of the matrix in (\ref{dPhi}) is:
\[
D\Phi (\mathbf c, \mathbf {\dot c}, \mu) = \left[
\begin{array}{ccc}
I & & \\
L_{21}
& I & \\
0 & h^T DF_{\mu}( \bc )^{-1} & 1 \\
\end{array}
\right]\
\left[
\begin{array}{ccc}
DF_{\mu}( \bc )
&0
&\dot F_{\mu}( \bc )
\\
&
DF_{\mu}( \bc )
&
U_{23}
\\
& &
U_{33}
\\
\end{array}
\right]
\]
where using (\ref{ift}):
{
\small
\begin{eqnarray*}
L_{21} &=&
\left(D^2F_{\mu}( \bc ) \otimes \bcdot\right)
DF_{\mu}( \bc )^{-1}
+
D\dot F_{\mu} (\bc) DF_{\mu}( \bc )^{-1}
\\
U_{23} &=&
2 D\dot F_{\mu}( \bc ) \ \bcdot
+
\ddot F_{\mu}( \bc )
+
D^2F_{\mu}( \bc ) (\bcdot, \bcdot)
\\
U_{33} &=&
-h^T \left(
2 DF_{\mu}( \bc )^{-1}
D\dot F_{\mu}( \bc ) \ \bcdot
+
DF_{\mu}( \bc )^{-1}
\ddot F_{\mu}( \bc )
+
DF_{\mu}( \bc )^{-1}
D^2F_{\mu}( \bc ) (\bcdot, \bcdot)
\right).
\end{eqnarray*}
}
Note that by construction, $F_{\mu} ({\mathbf c}(\mu)) \equiv 0$.
Differentiating once with respect to $\mu$, we obtain (\ref{ift}-2).
Differentiating once again,
\[
DF_{\mu} ({\mathbf c}(\mu))\  \mathbf{\ddot c}
+D^2 F_{\mu} ({\mathbf c}(\mu))\ (\mathbf{\dot c},\mathbf{\dot c})
+2 D\dot F_{\mu} ({\mathbf c}(\mu))\  \mathbf{\dot c} +
\ddot F_{\mu} ({\mathbf c}(\mu)) = 0
\]
Solving for $\bcdotdot$ and replacing into $U_{33}$, we
obtain:
\[
U_{33} = h^T \bcdotdot \ .
\]
We need to show that $U_{33} \ne 0$. Our hypothesis was
that $\dot \ga (t) \not \in \CH$. Multiplying equation
(\ref{curv-t}) by $h^T$ to the left, we obtain:
\[
h^T \dot \ga (t) = \frac{ h^T \bcdotdot(t)}
{\| \dot {\bc}(t) \|^2}.
\]
Hence, $U_{33}$ does not vanish and $D\Phi$ is non-singular
at $(\bc_0, \bcdot_0, \mu_0)$.
\end{proof}

\section{A B\'ezout bound for multi-homogeneous systems.}

According to Theorem \ref{th-length} to estimate the length of a
curve we have to count the number of points in a certain set. To
give such an estimate we use  the multi-homogeneous B\'ezout
Theorem. While this theorem is well-known to algebraic geometers,
topologists and homotopy method theorists, the computation of the
B\'ezout number is usually only carried out in the bi-homogeneous
case in textbooks. Morgan and Sommese \cite{mor1} prove the
theorem and give a simple description of how to compute the
number, which we repeat here.

Let $f = (f_i)_{1 \leq i \leq n}$ be a system of $n$ complex
polynomial equations in $n+m$ complex variables. These variables
are partitioned into $m$ groups $X_1, \ldots ,X_m$ with $k_j + 1$
variables into the $j-$th group. $f_i$ is said multi-homogeneous
if for any index $j$ there exists a degree $d_{ij}$ such that, for
any scalar $\la \in \C$,
$$f_i(X_1, \ldots , \la X_j, \ldots , X_m) = \la^{d_{ij}} f_i(X_1,
\ldots , X_j, \ldots , X_m).$$ In this case the system $f$ is
called multi-homogeneous. The B\'ezout number $\CB$ associated
with this system and this structure is defined as the coefficient
of $\prod_{j=1}^m \zeta_j^{k_j}$ in the product $\prod_{i=1}^n
\sum_{j=1}^m d_{ij}\zeta_j.$

We say that $(X_1, \ldots , X_m) \in \C^{n+m}$ is a zero for $f$
when $f(X_1, \ldots , X_m) = 0$. In that case, $f(\la_1 X_1,
\ldots , \la_m X_m) = 0$ for any $m-$tuple of complex scalars
$(\la_1, \ldots , \la_m)$. For this reason it is convenient to
associate a zero to a point in the product of projective spaces
$\P^{k_1}(\C) \times \ldots \times \P^{k_m}(\C)$. We use the same
notation for a point in $\P^{k_1}(\C) \times \ldots \times
\P^{k_m}(\C)$ and for any representative $(X_1, \ldots , X_m) \in
\C^{n+m}$.

We say that a zero $(X_1, \ldots , X_m) \in \P^{k_1}(\C) \times
\ldots \times \P^{k_m}(\C)$ is non-singular when the derivative
$$Df(X_1, \ldots , X_m) : \C^{n+m} \ra \C^{n}$$
is surjective. Notice that this definition is independent of the
representative $(X_1, \ldots , X_m) \in \C^{n+m}$. We have

\begin{theorem} \label{th-bezout} (Multi-homogeneous B\'ezout
Theorem) Let $f$ be a multi-ho\-mo\-ge\-ne\-ous system. Then the number of
isolated zeros of $f$ in $\P^{k_1}(\C) \times \ldots \times
\P^{k_m}(\C)$ is less than or equal to $\CB$. If all the zeros are
non-singular then $f$ has exactly $\CB$ zeros.
\end{theorem}

\section{The total curvature of the central path on the average.}
\label{formal}

To the matrix $A$ and the vector $b$ we {associate} the set of
admissible points of the primal problem via the set of
equalities-inequalities
$$Ax-s=b, \hspace{2em}
s\ge 0.$$ We may also consider the other polyhedra contained in
the subspace $Ax-s=b$ and defined by the inequalities
\[
s_i \ \ep_i \
0, \hspace{2em}1 \le i \le m,
\]
where $\ep = (\ep_1 , \ldots , \ep_n)$ is one of the $2^m$ vectors of sign conditions.

Let $\CF(A,b)$ denote the set of such primal strictly feasible
polyhedra contained in the subspace $Ax-s = b$  and $\CQ(A,b)$ the
set of those which are compact.

\begin{lemma}\label{Q} For $A$ and $b$ {almost everywhere},
\[
\# \CQ(A,b) = R_K(m,n) = \binomial{m-1}{n} \ .
\]
\end{lemma}
\proof {This statement was proved by Buck~\cite{buc} for $A$ and $b$
in {\em general position}. In particular, $A$ and $b$ are in general position
except in a set of measure zero, Lemma~\ref{Q} holds for $A$ and $b$ almost
everywhere.} \qed \medskip

\begin{proposition} A probability measure on \emph{P} which is
absolutely continuous with respect to Lebesgue measure is full.
\end{proposition}

\proof The set of $(A,b,c)$ in $\emph{P}$ where $(A,b)$ is not in
general position has zero Lebesgue measure by the above lemma and
Fubini's theorem, thus it has zero measure for any measure
absolutely continuous with respect to Lebesgue.\qed\medskip

Now we prove the corollary of the introduction assuming the Main
Theorem.

\proof The group $\emph{D}$ acts freely on $\emph{P}$, let
$\frac{\emph{P}}{\sim}$ denote the orbit space. Then we may
decompose the measure $\mu$ on the orbits of $\emph{D}$. Since
$\mu$ is sign invariant each point in the orbit gets equal measure
and the same is true for the conditional measure $\nu$, ie each
strictly feasible polytope in the orbit of $\emph{D}$ gets
equal measure when the measure $\nu$ is decomposed on orbits. Now
we average over the orbits of points in general position, apply
the Main Theorem and then average over $\frac{\emph{P}}{\sim}.$
\qed \medskip


It remains to prove Theorem \ref{th-main-sum}.

The proof of this theorem requires Lemmas~\ref{lem-main-th-1},
\ref{lem-main-th-2} and Proposition~\ref{prop-mu} below.

\begin{lemma} \label{lem-main-th-1} For each $\CF \in \CF (A,b)$, with non-empty primal central paths
 the Gauss curves associated with the central paths $\bc_{\primaldual}(\CF),$
$\bc_\primal(\CF)$and $\bc_\dual(\CF)$ are well-defined.
\end{lemma}

\begin{proof} The primal/dual (resp. primal, resp. dual) central path associated with a polyhedron $\CF \in \CF(A,b)$
satisfies the system of polynomial equations
\begin{equation}
F_{\mu}(x,s,y) =\left[
\begin{array}{l}
Ax-s - b \cr
A^Ty - c \cr
sy - \mu \e \cr
\end{array}
\right]
=0
\label{E1}
\end{equation}
with $\mu > 0$, and this system is the same for all those
polyhedra.

Let $D_s$ denote the diagonal matrix with diagonal entries $s_i$. Since $sy=\mu e$ (equation
\ref{E1}-3), $D_s$ is invertible. The derivative of $F_{\mu}$ is equal to
\[
DF_{\mu}(x,s,y) = \left[
\begin{array}{ccc}
A & -I & 0\cr
0  & 0   & A^T\cr
0  & D_y  & D_s\cr
\end{array}
\right]
\]
and it factors as:
\[
\left[
\begin{matrix}
I & 0& 0\\
0  &0 &I \\
0  &I&0 \\
\end{matrix}
\right]
\left[
\begin{array}{ccc}
I & & \cr
-\mu D_s^{-1} & I &\cr
0  &A^T D_s^{-1} & I\cr
\end{array}
\right]
\left[
\begin{array}{ccc}
-I &0 & A\cr
  & D_s & \mu D_s^{-1} A \cr
    & & -\mu A^T D_s^{-2} A \cr
\end{array}
\right]
\left[
\begin{matrix}
0 &I&0\\
0 &0 &I \\
I&0 & 0\\
\end{matrix}
\right].
\]
Therefore, since $A$ has full column rank, $\mu > 0$ and $s_j \ne 0$, this derivative is nonsingular
and we are in the conditions of the
Implicit Function Theorem. The speed vector
\[
\bcdot = (\dot x, \dot s, \dot y)
= -DF_{\mu}(x(\mu),y(\mu),z(\mu))^{-1}
\dot F_{\mu}((x(\mu),y(\mu),z(\mu))
\]
is the unique solution of the implicit equations:

\begin{equation}
\left\{
\begin{array}{l}
A\dot x - \dot s= 0 \cr
A^T \dot y = 0 \cr
\dot s y + s \dot y = \e \ .\cr
\end{array}
\label{E2}
\right.
\end{equation}
\medskip
\par
The Gauss curve for the primal-dual problem is
${(\dot x,\dot s,\dot y)}/{\| (\dot x,\dot s,\dot y) \|}$.
Notice that because of (\ref{E2}-3), $\dot s$, and $\dot y$ cannot be
together equal to $0$
so that this curve is well-defined.

The Gauss curve associated to the primal (resp. dual) central
path is
${(\dot x,\dot s)}/{\| (\dot x,\dot s) \|}$ (resp.
${\dot y}/{\| \dot y \|}$). Those curves are well defined,
for suppose that $\dot s = 0$. Then equations (\ref{E1}-3)
and (\ref{E2}-3) combined give:
\[
s y=\mu s \dot y.
\]
Hence, dividing componentwise by $s$ and then
multiplying by $A^T$, one obtains:
\[
c = A^T y = \mu A^T \dot y = 0
\]
which contradict the hypothesis $c \ne 0$. Suppose now $\dot y 0$. Then, by the same reasoning one obtains $s = \mu \dot s$.
Hence, by (\ref{E2}-1), $s$ is in the image of $A$. Then by
(\ref{E1}-1), $b$ is in the image of $A$, and hence the polyhedron
$Ax - s = b,$ $s \ge 0$ is either one point or unbounded. Thus, we
showed that the Gauss curves for the primal-dual, primal and dual
central paths are well-defined.
\end{proof}

{

A point of the curve $\gamma_{\primaldual}$ is the image under the
map
\[
(x,s,y,\dot x, \dot s, \dot y) \mapsto \frac{(\dot x, \dot s, \dot
y)} {\|(\dot x, \dot s, \dot y)\|}
\]
of a point $(x,s,y,\dot x, \dot s, \dot y)$ satisfying  the
systems (\ref{E1}) and (\ref{E2}) for some $\mu > 0$. Similarly,
a point of the curve $\gamma_{\primal}$ (resp. $\gamma_{\dual}$)
is the image of such a point under the map
\[
(x,s,y,\dot x, \dot s, \dot y) \mapsto \frac{(\dot x, \dot s)}
{\|(\dot x, \dot s)\|}
\]
\[
\text{(resp.\ } (x,s,y,\dot x, \dot s, \dot y) \mapsto \frac{\dot
y}{\|\dot y\|}\text{).}
\]
}

%
%
%
%

\begin{lemma} \label{lem-main-th-2} 
  It is assumed as above that $\CF \in \CF (A,b)$ and that $\bc_*$ and
$\gamma_*$ are defined as above.
  Let $u \in \mathbb R^n$, $v \in \mathbb R^m$, $w \in \mathbb R^m$ be not
all zero.
\begin{enumerate}
\item Each transversal intersection of the Gauss curve $\gamma_{\primaldual}$
with the hyperplane
\[
\CH_{\primaldual} = \left\{ (\dot x, \dot s, \dot y) \ :
u^T \dot x + v^T \dot s + w^T \dot y = 0 \right\}
\]
is the
image of a nonsingular solution of the polynomial system
\begin{equation}
\Phi^{A,b,c,u,v,w} (x,s,y,\dot x, \dot s, \dot y, \mu) = \left[
\begin{array}{l}
Ax - s - b \cr
A^Ty - c\cr
sy - \mu \e\cr
A\dot x - \dot s \cr
A^T \dot y  \cr
\dot s y + s \dot y - \e\cr
u^T \dot x + v^T \dot s + w^T \dot y
\end{array}
\right] = 0
\label{E3}
\end{equation}
such that $\mu > 0$.
\item Let $w=0$. Each transversal intersection of the Gauss curve $\gamma_{\primal}$
with the hyperplane
\[
\CH_{\primal} = \left\{ (\dot x, \dot s) \ :
u^T \dot x + v^T \dot s = 0 \right\}
\]
is the image
of a nonsingular solution of the polynomial
system (\ref{E3}).
\item Let $u=0$ and $v=0$.
Each transversal intersection of the Gauss curve $\gamma_{\dual}$
with the hyperplane
\[
\CH_{\dual} = \left\{ \dot y\   :
w^T \dot y = 0 \right\}
\]
is the image
of a nonsingular solution of the polynomial
system (\ref{E3}).
\end{enumerate}
\end{lemma}

\begin{proof}

Part 1 is Lemma~\ref{transverse}, where $F_{\mu}$ and $DF_{\mu}$ are
computed in (\ref{E1}) and (\ref{E2}).

Part 2 follows from the fact that any transversal intersection of
$\gamma_{\primal}$ with the hyperplane $u^T \dot x + v^T \dot s = 0$
corresponds to a transversal intersection of $\gamma_{\primaldual}$ with
the hyperplane $u^T \dot x + v^T \dot s + 0^T \dot y = 0$. Indeed,
if $\gamma_{\primal}(\mu)= \frac{(\dot x, \dot s)}{\| \dot x , \dot s\|}$ we set
$\gamma_{\primaldual}(\mu) = \frac{(\dot x, \dot s, \dot y)}{\| \dot x , \dot s, \dot y \|}$. Then
$(u,v)^T \gamma_{\primal}(\mu) = 0$ if and only if
$(u,v,0)^T \gamma_{\primaldual}(\mu) = 0$.

Now, assume that the intersection of $\gamma_{\primal}$ with $u^T
\dot x + v^T \dot s = 0$ is transversal. Then,

\begin{eqnarray*}
\frac{\partial}{\partial \mu}
(u,v,0)^T \gamma_{\primaldual}(\mu)
&=&
\frac{1}{\| \dot x , \dot s, \dot y \|}
\left(u^T \ddot x + v^T \ddot s + 0^T \ddot y \right)
+
\left(
u^T \dot x + v^T \dot s + 0^T \dot y
\right)
\frac{\partial}{\partial \mu}
\frac{1}{\| \dot x , \dot s, \dot y \|}
\\
&=&
\frac{1}{\| \dot x , \dot s, \dot y \|}
\left(u^T \ddot x + v^T \ddot s + 0^T \ddot y\right)
\\
&=&
\frac{\| \dot x, \dot s \|}{\| \dot x , \dot s, \dot y \|}
\frac{\partial}{\partial \mu}
(u,v)^T \gamma_{\primal}(\mu)
\ne 0
\end{eqnarray*}
and therefore the intersection of $\gamma_{\primaldual}$ with $u^T
\dot x + v^T \dot s + 0^T \dot y = 0$ is also transversal.

The proof of Part 3 is similar.
\end{proof}

%
%

\begin{proposition} \label{prop-mu}  Let $m > n \ge 1$. Let $A$ be an $m \times n$
matrix of rank $n$, and let $b \in \R^m$ and $c \in \R^n$, $c$
non-zero. Then, for any transversal hyperplane $\CH_*$, the
polynomial system (\ref{E3}) has at most
\[
\CB_{\primaldual} \le 2 n \binomial{m-1}{n}
\]
nonsingular solutions $(x,s,y,\dot x, \dot s, \dot y, \mu) \in
\R^n \times \R^m \times \R^m \times \R^n \times \R^m \times \R^n
\times \R$ with $\mu > 0$.

If
furthermore we have $w=0$,
the number of nonsingular solutions is bounded
above by
\[
\CB_{\primal} \le 2 (n-1) \binomial{m-1}{n}
\]

If instead we have $u=0$ and $v=0$,
the number of nonsingular solutions is still bounded
above by
\[
\CB_{\dual} \le 2 n \binomial{m-1}{n}
\]
\end{proposition}

The proof of Proposition~\ref{prop-mu} is long, and is postponed
to Section~\ref{rootcount}.

%
%
%
%
%
%

\begin{proof}[Proof of Theorem~\ref{th-main-sum}]
The total curvature is the sum of the lengths of the Gauss curves
corresponding to strictly feasible regions . According to
Corollary \ref{cor-length}, a bound $\CB_{*}$ for the number of
intersections (counted with multiplicity) of the associated Gauss
curves with a transversal hyperplane gives the bound $\pi \CB_{*}$
for the length. Finally, by lemma \ref{lem-main-th-2} and
proposition \ref{prop-mu} $\CB_{*}$ may be taken as is proposition
\ref{prop-mu}.
\end{proof}

\section{Proof of Proposition \ref{prop-mu}}
\label{rootcount}

{ The proof of Proposition \ref{prop-mu} is quite long,
and occupies all of this section. There are actually three cases,
that are quite similar and will be treated in parallel. The symbol
$*$ stands for $\primaldual$, $\primal$ or $\dual$. Each of these
cases will be known as the {\em primal-dual}, the {\em primal} and
the {\em dual} case, respectively.

We proceed as follows:

\subsection{Complexification of the equations}
The first step is to complexify the equations,
i.e. to keep the coefficients fixed and to
consider the variables as complex instead of real.

\begin{lemma}

  The number of nonsingular solutions of (\ref{E3}) in $\mathbb R^{4m+2m+1}$
with $\mu > 0$ is bounded above by the number of nonsingular solutions of
(\ref{E3}) in $\mathbb C^{4m+2m+1}$ with $\mu \ne 0$.
\end{lemma}
\begin{proof}
 A real root is, in particular, a complex root. It
is non-degenerate if and only if the determinant of the Jacobian
matrix  of the derivative does not vanish. The non-vanishing of
this determinant does not depend on whether the matrix is
considered as real or complex.
\end{proof}

Note that when we complexify the equations, the terms $u^T \dot x +
v^T \dot s + w^T \dot y$ stand for the usual transpose.

A standard application of B\'ezout's Theorem implies that:
\begin{lemma} \label{finiteness} 
  The number of nonsingular solutions of (\ref{E3}) in
in $\mathbb C^{4m+2m+1}$ with $\mu \ne 0$ is bounded above
by $2^{2m}$.
\end{lemma}
This estimate, while
ensuring finiteness, is not sharp enough for our theorem.

\subsection{Continuation of non-degenerate roots}

More formally, we denote by $\mathcal A_{\primaldual}$ the set of
all complex $A,b,c,u,v,w$ where $A$ has rank $n$, $c \ne 0$ and
$u$, $v$, $w$ are not simultaneously zero. We also denote by
$\mathcal A_{\primal}$ (resp. $\mathcal A_{\dual}$) the
intersection of $\mathcal A_{\primaldual}$ with the linear space
$w=0$ (resp. $u=0$ and $v=0$).

Then, $\CB_*$ will denote the maximal number of non-degenerate complex
roots of (\ref{E3}) with $\mu \ne 0$, where $*$ is one of $\primaldual$,
$\primal$, $\dual$ and the maximum is taken over all parameters in $\mathcal A_*$.
As in Remark~\ref{finiteness}, $\CB_*$ is finite. Hence this maximal
number is attained, and at that point all the non-degenerate complex
roots may be continued in a certain neighborhood. Thus,

\begin{lemma}\label{open}
  The maximal number $\CB_*$ of non-degenerate complex roots is attained
in a certain open set of $\mathcal A_*$.
\end{lemma}

\begin{proof}

Lemma~\ref{finiteness} implies that $\CB_*$ is attained for some
parameter $(A,b,c,u,v,w)$.

By the Implicit Function Theorem, the non-degenerate complex roots
of $(\ref{E3})$ with $\mu \ne 0$ can be continued to non-degenerate
complex roots with $\mu \ne 0$, in a certain neighborhood of the
parameter $(A,b,c,u,v,w)$.
\end{proof}

\subsection{Non-degeneracy at the maximum}

The following fact will be needed in the sequel:

\begin{proposition}\label{non-degeneracy}
  The complex roots of (\ref{E3}) with $\mu \ne 0$ are all non-degenerate,
almost everywhere in $\mathcal A_*$.
\end{proposition}


\begin{proof}
Let $\mathcal X = \{ x, s, y,
\dot x, \dot s, \dot y, \mu \in \C^{4m+2n+1}: \mu \ne 0 \}$. We consider the
evaluation function
\[
\begin{array}{lrcl}
\mathrm{ev}:& \mathcal A_* \times \mathcal X & \rightarrow & \mathbb C^{4m+2n+1} \\
& (A,b,c,u,v,w ; x,s,y, \dot x, \dot s, \dot y,\mu) &
\mapsto &
\Phi^{A,b,c,u,v,w} (x,s,y, \dot x, \dot s, \dot y, \mu)
\end{array}
\]
where $\Phi$ was defined in (\ref{E3}).

$0$ is a regular value of $\mathrm{ev}$ if and only if $D\mathrm{ev}(A,b,c,\dots, \mu)$
is onto when $\mathrm{ev}(A,b,c,\dots, \mu) = 0$ (See~\cite[Ch II \S 1]{gol}).

\begin{lemma}\label{regval}
$0$ is a regular value for $\mathrm{ev}$.
\end{lemma}

This
Lemma guarantees that $V = \mathrm{ev}^{-1}(0)$ is a smooth
manifold and $\dim V = \dim \mathcal A_*$.

Now we consider the natural projection $\pi_1: V \rightarrow
\mathcal A_*$. By Sard's theorem, the regular values of $\pi_1$
have full measure in $\mathcal A_*$. Since $\dim V = \dim A_*$,
$(A, b, \dots, w)$  is a regular value if and only if $D\pi_1$ is
an isomorphism at every  point $(A, b, \dots, \mu)$ such that
$\pi_1(A, b, \dots, \mu)=(A, b, \dots, w)$. For such systems, all
the roots with $\mu \ne 0$ are non-degenerate.
\end{proof}

\begin{proof}[Proof of Lemma~\ref{regval}]

We first reorder the equations and the variables of (\ref{E3}) as follows:

\begin{equation}
\Phi (b,c,u,v,w,y,\dot s, \dot y, \dot x, x, s, A, \mu) = \left[
\begin{array}{l}
Ax - s - b \cr
A^Ty - c\cr
u^T \dot x + v^T \dot s + w^T \dot y \cr
sy - \mu \e\cr
A\dot x - \dot s \cr
\dot s y + s \dot y - \e\cr
A^T \dot y  \cr
\end{array}
\right]
\label{E3bis}
\end{equation}

In order to show that $D\Phi$ has full rank $4m+2n+1$, we will show that a
certain submatrix has rank $4m+2n+1$. Namely, we will consider only
the derivatives with respect to variables $b$ to $\dot x$, and
derivation with respect to $x$, $s$, $A$ and $\mu$ will be omitted.
We obtain the block matrix
\[
D_{b, \cdots, \dot x}\Phi = \left[
\begin{matrix}
-I &    &             &            &     &     & \\
   & -I &             & A^T        &     &     & \\
   &    & \dot x^T \dot s^T \dot y^T
                      &            & v^T & w^T & u^T \\
   &    &             & D_s        &     &     & \\
   &    &             &            & -I  &     & A \\
   &    &             & D_{\dot s} & D_y & D_s & & \\
   &    &             &            &     & A^T & & \\
\end{matrix}
\right]
\]

Recall that $\mu \ne 0$, hence no coordinate of $s$ or $y$ can
vanish and the diagonal matrices $D_s$ and $D_y$ have full rank.

Performing row operations on the previous matrix, one obtains:
\[
D_{b, \cdots, \dot x}\Phi = L
\left[
\begin{matrix}
-I &    &             &            &     &     & \\
   & -I &             & A^T        &     &     & \\
   &    & \dot x^T \dot s^T \dot y^T
                      &            & v^T & w^T & u^T \\
   &    &             & D_s        &     &     & \\
   &    &             &            & -I  &     & A \\
   &    &             &            &     & D_s & D_y A \\
   &    &             &            &     &     & -A^T D_{s}^{-1}D_y A \\
\end{matrix}
\right]
\]
for an invertible lower triangular matrix $L$.
Since not all of $\dot x_i$, $\dot s_i$, $\dot y_i$ can be
zero (Lemma~\ref{lem-main-th-1}) and $D_s$ has full rank,
it remains only to check that $-A^T D_{s}^{-1}D_y A$ has
also full rank. This follows from the identity:
\[
-A^T D_{s}^{-1}D_y A = \mu ( D_{s}^{-1} A)^T ( D_{s}^{-1} A)
\]
and from the fact that $A$ has full rank.

Hence, $D\Phi$ has rank $4m+2n+1$, and we are done.
\end{proof}

}
%
%

\subsection{Genericity}
{

In this section we show that it is sufficient to bound the number
of non-degenerate roots of systems satisfying conditions 1 through
5 of proposition \ref{genericity}.

 Let $K \subset \{ 1, \cdots, m\}$.
 We define $S_K$ as the linear space of all $\{s \in \mathbb
C^m : \ s_k = 0 \ \forall \ k \in K\}$.

\begin{proposition}\label{genericity}
  Let $m > n$. There is a point $(A,b,c,u,v,w) \in \mathcal A_*$ such that:
\begin{enumerate}

\item \label{gen-max} The maximal number $\CB_*$ of non-degenerate
complex solutions of (\ref{E3}) with $\mu \ne 0$ is attained at
this point. \item $b\neq 0.$\item \label{gen-nondeg} All the
solutions at that point are non-degenerate. \item
\label{gen-primal} For any $K \subset \{1, \cdots, m \}$, the
linear space $S_K$ and the affine space $(\ker A^H)^{\perp} - b$
intersect if and only if $n-\#K \ge 0$. In that case, the
intersection has dimension $n -\#K$. \item \label{gen-dual} For
any $K \subset \{1, \cdots, m \}$, the linear space $S_{\{1,
\dots, m\} \setminus K}$ and the affine space $\{ \ga: A^T \ga c\}$ intersect if and only if $\#K - n \ge 0$. In that case, the
intersection has dimension $\#K - n.$

\end{enumerate}
\end{proposition}

\begin{proof}

  By Lemma~\ref{open}, item 1 holds on an open set $U \subset \mathcal A^{*}$.
Items 2,3 will fail only on zero measure set
(Proposition~\ref{non-degeneracy}). For items 4 and 5, notice that
with probability one, $\dim (\ker A^H)^{\perp} - b = n$ and $\dim
\{ \ga: A^T \ga = c\} = m-n$. On the other hand, $\codim S_K \#K$ and $\codim S_{\{1, \dots, m\} \setminus K} = m - \#K$. Thus
it is easy to see that systems points violating items 3 and 4 will
fail in a finite union of zero measure sets.

Hence, items 2 to 5 will hold on a subset of $\mathcal A^{*}$ of
full measure which has a non-empty intersection with the open set
of Lemma~\ref{open}.
\end{proof}

This result has the following consequence: to give a bound
for the number of non-degenerate solutions of the system (\ref{E3})
with $\mu \ne 0$, we can replace
the initial data by the data $(A,b,c,u,v,w)$ of
Proposition~\ref{genericity}.

Also, for convenience,
we will count the number of
{\em isolated} roots of the corresponding system, which is the
same.

}
\subsection{Simplification of the equations}

\begin{lemma} \label{lem-main-th-4} Set $\hat u = u + A^T v$. The polynomial systems (\ref{E3}) and
(\ref{eq-mu}) below have the same solutions with $\mu \not = 0$ so that, the isolated solutions of
(\ref{E3}) with $\mu \not = 0$ are identical to the isolated solutions of (\ref{eq-mu}) with $\mu \not = 0$.
\begin{equation}
\Psi_\mu (x,s,y,\dot x, \dot s, \dot y) = \left[
\begin{array}{l}
Ax - s - b \cr
A^Ty - c\cr
sy - \mu \e\cr
A\dot x - \dot s \cr
A^T \dot y  \cr
\mu \dot y + y^2 \dot s - y\cr
\hat u^T \dot x + w^T \dot y
\end{array}
\right] = 0.
\label{eq-mu}
\end{equation}
\end{lemma}

\begin{proof}
This last system is obtained from (\ref{E3}) by the transformation
\begin{equation}
\left\{
\begin{array}{l}
\Psi_{k} = \Phi_{k}, \ 1 \le k \le 5, \cr
\Psi_{6} = y(\Phi_{6} - \Phi_{3}),\cr
\Psi_{7} = \Phi_{7} + v^T \Phi_{4}. \cr
\end{array}
\right.
\end{equation}
When $\mu \not = 0$, then no component of $y$ is zero so the
solutions of (\ref{E3}) and the solutions of
(\ref{eq-mu})coincide.
\end{proof}

\medskip
Let
$g_1, \cdots, g_{m-n}$ and $f \in \C^m$, $\ga$, $\de$ and $\hat w \in \C^{m-n}$
be such that
\begin{enumerate}
\item[(a)] $(g_1, \cdots, g_{m-n})$ is a basis of $\ker A^T$,
\item[(b)] $A^T f = c$,
\item[(c)] $y = f + \sum_{j=1}^{m-n} \ga_j g_j,$
\item[(d)] $\dot y = \sum_{j=1}^{m-n} \de_j g_j,$
\item[(e)] $\hat w_j = w^T g_j$, $\ j=1 \ldots m-n.$
\end{enumerate}
We will denote by $g_{ij}$ the $i$-th coordinate of $g_j$, by $G$ the $m \times (m-n)$ matrix with entries $g_{ij}$, by $A_i$ the $i-$th row of the matrix $A$, $i=1, \dots, m,$ and by {$A^\da = (A^HA)^{-1}A^H$}
the Moore-Penrose inverse of $A$ (injective); $A A^\da$ is the orthogonal projection onto $\im A$ while $A^\da A = I_n$.

\begin{lemma} \label{lem-main-th-5} The system \ref{eq-mu} has the same solutions as
\begin{equation}
\Om_\mu (x,s,y,\dot x, \dot s, \dot y) = \left[
\begin{array}{l}
x - A^\da(s + b) \cr
g_j^T (b+s), \ j=1 \ldots m-n, \cr
A^Ty - c\cr
sy - \mu \e\cr
A\dot x - \dot s \cr
A^T \dot y  \cr
\mu \dot y + y^2 \dot s - y\cr
\hat u^T \dot x + w^T \dot y
\end{array}
\right] = 0
\label{eq-om}
\end{equation}
\end{lemma}

\begin{proof} The equations \ref{eq-om}-2 are equivalent to $b+s \in \left(\ker A^T \right)^\perp$ that is $b+s \in \im A$. Under this assumption, $x = A^\da(s + b)$ gives $Ax = b+s$. Thus \ref{eq-om}-1 and \ref{eq-om}-2 imply \ref{eq-mu}-1. Conversely, if $Ax = b+s$ we get $b+s \in \im A$ and $x = A^\da(s + b)$ that is
\ref{eq-om}-1 and \ref{eq-om}-2.
\end{proof}

\subsection{Elimination of variables}
Let us now introduce a new polynomial system with the same number of zeros:

\noindent $\Xi_\mu (s,\dot x, \ga, \de) = $
\begin{equation}
\left[
\begin{array}{l}
g_j ^T (b+s) \rangle \cr
s_i (f_i + \sum_{j=1}^{m-n} \ga_j g_{ij}) - \mu \cr
\mu \sum_{j=1}^{m-n} \de_jg_{ij} + \left(f_i + \sum_{j=1}^{m-n} \ga_j g_{ij}\right)^2 A_i \dot x - \left(f_i + \sum_{j=1}^{m-n} \ga_j g_{ij}\right) \cr
\hat u^T \dot x + \hat w^T \de \cr
\end{array}
\right] = 0,
\label{eq-xi}
\end{equation}
where the range for $j$ in \ref{eq-xi}-1 is $\{ 1, \ldots ,m-n \}$ and the range for $i$ in \ref{eq-xi}-2 and \ref{eq-xi}-3 is  $\{ 1, \ldots ,m \}$.
We also define the following maps:
$$\La : \C^{m} \times \C^{n} \times \C^{m-n} \times \C^{m-n} \times \C \ra \C^{n} \times \C^{m} \times \C^{m} \times \C^{n} \times \C^{m} \times \C^{m} \times \C$$
$$\La(s,\dot x, \ga, \de, \mu) = \left( A^\da(s + b), s, f + \sum_{j=1}^{m-n} \ga_j g_j, \dot x, A\dot x, \sum_{j=1}^{m-n} \de_j g_j, \mu
\right),$$
and the projection
$$\Pi_{2478} : \C^{n} \times \C^{m-n} \times \C^{n} \times \C^{m} \times \C^{m} \times \C^{n} \times \C^{m} \times \C \ra
\C^{m-n} \times \C^{m} \times \C^{m} \times \C$$
$$\Pi_{2478}(z_1,z_2,z_3,z_4,z_5,z_6,z_7,z_8) = (z_2,z_4,z_7,z_8).$$

\begin{lemma} \label{lem-main-th-6}
The construction of the system (\ref{eq-xi}) is such that the diagram
$$\begin{array}[pos]{rcccl}
&            &       \Om     &            & \cr
&\C^{4m+2n+1}&\longrightarrow&\C^{4m+2n+1}& \cr
\La&    \uparrow & & \downarrow & \Pi_{2478}\cr
     &\C^{3m-n+1}&{\longrightarrow}&\C^{3m-n+1}& \cr
     &            &       \Xi     &            & \cr
\end{array}$$
is commutative. Moreover
\begin{enumerate}
\item If $(s,\dot x, \ga, \de, \mu)$ is a solution of (\ref{eq-xi}) then
$(x,s,y,\dot x, \dot s, \dot y, \mu)$ is a solution of (\ref{eq-om}),
\item Any solution $(x,s,y,\dot x, \dot s, \dot y, \mu)$ of (\ref{eq-om})
is equal to $\La(s,\dot x, \ga, \de, \mu)$ for a unique solution of (\ref{eq-xi}),
\item Any isolated solution of (\ref{eq-om}) with $\mu \ne 0$ corresponds to an
isolated solution of (\ref{eq-xi}) with $\mu \ne 0$.
\end{enumerate}
\end{lemma}

\begin{proof} The proof is easy and left to the reader.
\end{proof}

\medskip
Now we look at equations \ref{eq-xi}-3, \ref{eq-xi}-4 as a linear system of $m+1$ equations in the $m$ unknowns $(\de, \dot x)$, with coefficients depending on $\ga$ and $\mu$. When $(s, \dot x, \ga, \de, \mu)$ is a solution of (\ref{eq-xi}) with $\mu \ne 0$, those equations have a solution if and only if the determinant of the corresponding augmented matrix vanishes. We can write the augmented matrix as:
{\small
\begin{equation} \label{eq-mu-aug}
M(\ga, \mu)= \left[
\begin{matrix}
\mu g_{11} & \cdots & \mu g_{1,m-n} & a_{11}Q_1(\ga) & \cdots & a_{1n} Q_1(\ga)& f_1+\sum_{j=1}^{m-n} \ga_j g_{1j} \\
\vdots &        & \vdots    & \vdots &        & \vdots  & \vdots \\
\vdots &        & \vdots    & \vdots &        & \vdots  & \vdots \\
\vdots &        & \vdots    & \vdots &        & \vdots  & \vdots \\
\mu g_{m1} & \cdots & \mu g_{m,m-n} & a_{m1}Q_m(\ga) & \cdots & a_{mn} Q_m(\ga)& f_m+\sum_{j=1}^{m-n} \ga_j g_{mj} \\
\hat w_1 & \cdots & \hat w_{m-n} &  \hat u_1   & \cdots & \hat u_n     & 0 \\
\end{matrix}
\right]
\end{equation}}
where $Q_{i}(\ga) = \left(f_i + \sum_{j=1}^{m-n} \ga_j g_{ij}\right)^2$.
\medskip
By dividing the first $m-n$ columns by $\mu$ and then multiplying the
last row by $\mu$, we obtain a new matrix $M'$ whose determinant
vanishes if and only if the determinant of $M$ vanishes:

{\small
\begin{equation}\label{eq-mu-aug3}
M'(\ga, \mu) = \left[
\begin{matrix}
g_{11} & \cdots & g_{1,m-n} & a_{11}Q_1(\ga) & \cdots & a_{1n} Q_1(\ga)& f_1+\sum_{j=1}^{m-n} \ga_j g_{1j} \\
\vdots &        & \vdots    & \vdots &        & \vdots  & \vdots \\
\vdots &        & \vdots    & \vdots &        & \vdots  & \vdots \\
\vdots &        & \vdots    & \vdots &        & \vdots  & \vdots \\
g_{m1} & \cdots & g_{m,m-n} & a_{m1}Q_m(\ga) & \cdots & a_{mn} Q_m(\ga)& f_m+\sum_{j=1}^{m-n} \ga_j g_{mj} \\
\hat w_1 & \cdots & \hat w_{m-n} &  \mu \hat u_1   & \cdots & \mu \hat u_n     & 0\\
\end{matrix}
\right]
\end{equation}}
Its determinant is the same as the determinant of:

{
\small
\begin{equation} \label{eq-mu-aug2}
M''(\ga, \mu) = \left[
\begin{matrix}
g_{11} & \cdots & g_{1,m-n} & a_{11}Q_1(\ga) & \cdots & a_{1n} Q_1(\ga)& f_1 \\
\vdots &        & \vdots    & \vdots &        & \vdots  & \vdots \\
\vdots &        & \vdots    & \vdots &        & \vdots  & \vdots \\
\vdots &        & \vdots    & \vdots &        & \vdots  & \vdots \\
g_{m1} & \cdots & g_{m,m-n} & a_{m1}Q_m(\ga) & \cdots & a_{mn} Q_m(\ga)& f_m \\
\hat w_1 & \cdots & \hat w_{m-n} &  \mu \hat u_1   & \cdots & \mu \hat u_n & - \sum \ga_i \hat w_i\\
\end{matrix}
\right]
\end{equation}}

{ We now have to distinguish the three cases of
Proposition~\ref{prop-mu}. In the primal-dual case, we define the
eliminating
polynomial $h_{\primaldual}(\ga,\mu) =  \det M''(\ga, \mu)$.
In the dual case, we also define
$h_{\dual}(\ga) =  \det M''(\ga, \mu)$ but now, since $\hat u=0$,
the eliminating polynomial is independent of $\mu$.
In the primal case, $\hat w = 0$ and we notice that the last row of
$M''$ is divisible by $\mu$. Hence, we set $h_{\primal}(\ga)= \mu^{-1} \det M''(\ga, \mu)$.

\begin{lemma} \label{lem-main-th-7}  With the notations above,
$(s, \dot x, \ga, \de, \mu)$ is a solution of (\ref{eq-xi}) with $\mu \ne 0$ if and only if
\begin{equation}\Upsilon (s, \dot x, \ga, \de, \mu) = \left[
\begin{array}{l}
s_1 (f_1 + \sum_{j=1}^{m-n} \ga_j g_{1j}) - \mu \cr
M(\gamma, \mu) \left(
\begin{array}{c}
    \delta \cr
    \dot x \cr
    1 \cr
\end{array}
\right)\cr
G^T(b+s)  \cr
s_i (f_i + \sum_{j=1}^{m-n} \ga_j g_{ij}) -  s_1 (f_1 + \sum_{j=1}^{m-n} \ga_j g_{1j})\cr
h_*(\gamma,s_1)\cr
\end{array}
\right] = 0,
\label{eq-Up}
\end{equation}
and $\mu \ne 0$.
\end{lemma}

\begin{proof} The system $\Upsilon = 0$ is the same as $\Xi = 0$ plus the equation $h(\ga , s_1) = 0$ which is a
consequence of $\Xi = 0$ as has been explained previously.
\end{proof}

\begin{lemma} \label{lem-main-th-8} The number of isolated solutions of the system (\ref{eq-Up}) with $\mu \ne 0$ is less than or equal to the number of isolated solutions with $s_1 (f_1 + \sum_{j=1}^{m-n} \ga_j g_{1j}) \ne 0$ of

\begin{equation} \Theta (s, \ga) = \left[
\begin{array}{l}
G^T(b+s)  \cr
s_i (f_i + \sum_{j=1}^{m-n} \ga_j g_{ij}) -  s_1 (f_1 + \sum_{j=1}^{m-n} \ga_j g_{1j})\cr
h_*(\gamma,s_1)\cr
\end{array}
\right] = 0
\label{eq-Th}
\end{equation}
where the range for $i$ in the second equation is $2 \le i \le m$.
\end{lemma}

\begin{proof} This lemma is a consequence of the following facts:
\begin{itemize}
\item An isolated solution $(s, \dot x, \ga, \de, \mu)$ of (\ref{eq-Up}) gives an isolated solution
$(s, \ga)$ of (\ref{eq-Th}),
\item Two distinct solutions $(s, \dot x, \ga, \de, \mu)$ and $(s', \dot x', \ga', \de', \mu')$ of (\ref{eq-Up}) with $\mu \ne 0$ and $\mu' \ne 0$ give two distinct solutions $(s, \ga)$ and $(s', \ga')$ of (\ref{eq-Th}) with
$s_1 (f_1 + \sum_{j=1}^{m-n} \ga_j g_{1j}) \ne 0$ and $s_1' (f_1 + \sum_{j=1}^{m-n} \ga_j' g_{1j}) \ne 0$.
\end{itemize}

The first fact is true because (\ref{eq-Th}) is a sub-system of (\ref{eq-Up}).

Let us prove the second assertion. Let $(s, \dot x, \ga, \de, \mu)$ and $(s, \dot x', \ga, \de', \mu')$ be two solutions of (\ref{eq-Up}) with $s_1 (f_1 + \sum_{j=1}^{m-n} \ga_j g_{1j}) \ne 0$. Our aim is two prove that
$(\dot x, \de, \mu) = (\dot x', \de', \mu')$. We have clearly $\mu = \mu' \ne 0$ and $s_i (f_i + \sum_{j=1}^{m-n} \ga_j g_{ij}) \ne 0$ for each $i$. Moreover $(\dot x, \de)$ is given by the system
$$ M(\ga,\mu)
\left(
\begin{array}[pos]{c}
    \de\\
    \dot x\\
    -1\\
\end{array}
\right)
\left(
\begin{array}[pos]{ccc}
    \mu G&\diag (Q_i(\ga))A&f+G\ga\\
    \hat w&\hat u&0\\
\end{array}
\right)
\left(
\begin{array}[pos]{c}
    \de\\
    \dot x\\
    -1\\
\end{array}
\right) = 0.
$$
This system has a unique solution if and only if $\rank M(\ga,\mu) = m$. Let us prove that the $m \times m$ submatrix
$$\left(
\begin{array}[pos]{cc}
    \mu G&\diag (Q_i(\ga))A\\
\end{array}
\right)$$
is nonsingular. Let $\al \in \C^{m-n}$ and $\be \in \C^n$ be given. If
$$\left(
\begin{array}[pos]{cc}
    \mu G&\diag (Q_i(\ga))A\\
\end{array}
\right)
\left(
\begin{array}[pos]{c}
    \al\\
    \be\\
\end{array}
\right) = 0
$$
then
$$\mu G \al + \diag (Q_i(\ga))A \be = 0$$
so that, multiplying by $A^T$, we get
$$ A^T \diag (Q_i(\ga))A \be = 0.$$
Since $(f_i + \sum_{j=1}^{m-n} \ga_j g_{ij})^2 > 0$ and $\rank A = n$ the matrix $ A^T \diag (Q_i(\ga))A $ is positive definite and
consequently $\be = 0$. Thus $\mu G \al = 0$ so that $\al = 0$
because $\mu \ne 0$ and $\rank G = m-n$. This completes the proof.
\end{proof}
}

\subsection{The B\'ezout bound}
To count the number of isolated roots of this last system we use the multi-homogeneous B\'ezout bound given in Theorem \ref{th-bezout}. For this purpose we bi-homogenize the system (\ref{eq-Th}): we introduce the homogenizing variables $s_0$ in the first group of variables and $\ga_0$ in the second group. We obtain the system:
\begin{equation} \Theta (s_0,s,\ga_0,\ga) = \left[
\begin{array}{l}
G^T(s_0b+s)  \cr
s_i (\ga_0 f_i + \sum_{j=1}^{m-n} \ga_j g_{ij}) - s_1 (\ga_0 f_1 + \sum_{j=1}^{m-n} \ga_j g_{1j}) =0 \cr
h_*(s_0,s,\ga_0,\ga) = 0 \ .
\end{array}
\right] = 0.
\label{eq-Thh}
\end{equation}

{
Here $h_*(s_0,s,\ga_0,\ga)$ is just the homogeneization of $h_*(s,\ga)$.
For instance, $h_{\primaldual}$ is
the determinant of the matrix obtained from $M''(\ga, \mu)$ by the following changes:
\begin{enumerate}
\item $Q_i(\ga)$ becomes $Q_i(\ga_0, \ga) = (\ga_0 f_i + \sum_{j=1}^{m-n} \ga_j g_{ij})^2$,
\item $\mu$ becomes $\ga_0 f_1 + \sum_{j=1}^{m-n} \ga_j g_{1j}$,
\item $f_i$ in the last column becomes $\ga_0 f_i$,
\item $- \sum \ga_i \hat w_i$ in the last column becomes $- s_0 \sum \ga_i \hat w_i$.
\end{enumerate}

\begin{lemma} \label{lem-main-th-9} \ \begin{enumerate}
\item The system \ref{eq-Thh} is bi-homogeneous in the groups of variables
$(s_0,s) \in \C^{m+1}$ and $(\ga_0,\ga) \in \C^{m-n+1}$,
\item There are $m-n$ equations of multi-degree $(1,0)$ and $m-1$ equations of
multi-degree $(1,1)$. The multi-degree of the last equation depends on whether
$*$ stands for $\primaldual$, $\primal$ or $\dual$:
\subitem The multi-degree of $h_{\primaldual}$ is $(1,2n+1)$.
\subitem The multi-degree of $h_{\primal}$ is $(0,2n-2)$.
\subitem The multi-degree of $h_{\dual}$ is $(0,2n+1)$.
\item $(s,\ga)$ is an isolated solution of (\ref{eq-Th}) if and only if $(1,s,1,\ga)$ is an isolated solution of (\ref{eq-Thh}),
\item The number of isolated solutions of
(\ref{eq-Thh}) is bounded as follows:
\subitem When $*$ is equal to $\primaldual$, by
$2 n \binomial{m-1}{n} + \binomial{m}{n}$,
\subitem when $*$ is equal to $\primal$, by
$(2n-2) \binomial{m-1}{n}$ ,and
\subitem when $*$ is equal to $\dual$, by
$(2n+1) \binomial{m-1}{n}$.
\end{enumerate}
\end{lemma}

\begin{proof}  The assertions 1 and 2 hold by construction.

The third assertion comes from a classical fact: when both product spaces are equipped with the usual topology, the canonical injection
$$i : \C^m \times \C^{m-n} \ra \P^{m}(\C) \times \P^{m-n}(\C)$$
is open and continuous.

By the Multi-homogeneous B\'ezout Theorem (Theorem~\ref{th-bezout}), the number of isolated solutions in the primal-dual case ($* = \primaldual$)
is no more than the coefficient of
$\zeta_1^{m} \zeta_2^{m-n}$ in the expression
\[
\zeta_1^{m-n}\ (\zeta_1 + \zeta_2)^{m-1}\ \left(\zeta_1+(2n+1)\ \zeta_2 \right).
\]
In the primal-dual case (item 6), this coefficient is precisely
\[
(2 n+1) \binomial{m-1}{n} + \binomial{m-1}{n-1}
2 n \binomial{m-1}{n} + \binomial{m}{n} \ .
\]
The primal and dual cases are similar.
\end{proof}
}
\subsection{The spurious roots}

Except for the primal case, the number of roots we have obtained
is still too big to give the bound announced in the Theorem
\ref{th-main-sum}. We have to eliminate some of them to obtain the
right result.

There are three classes of bi-homogeneous solutions $(s_0,s,\ga_0,\ga)$
of the system (\ref{eq-Thh}).
\begin{enumerate}
\item{\bf Roots at infinity} are roots for which $\ga_0=0$ or
$s_0 = 0$. We will not worry about roots at infinity here,
\item{\bf Spurious roots}. These are the {finite}
roots for which $\mu = 0$ that is
$s_1 = 0$ or $\ga_0 f_1 + \sum_{j=1}^{m-n} \ga_j g_{1j} = 0$,
\item{\bf Legitimate roots} are all the other solutions.
\end{enumerate}
Notice that solutions $(s,\ga, \mu)$ of the system (\ref{eq-Th}) with $s_1(f_1 + \sum_{j=1}^{m-n} \ga_j g_{1j}) \ne 0$
correspond always to legitimate roots of (\ref{eq-Thh}).

{

\begin{lemma}  \label{lem-main-th-10} The number
of spurious roots of
(\ref{E3}) for $*$ one of $\primaldual$ or $\dual$ is $\binomial{m}{n}$.
\end{lemma}

\begin{proof} 
We will only deal with the primal-dual case, the dual case being
exactly the same.

The idea of the proof is to produce a bijection
from the spurious roots to the class of
subsets $K \subset \{1, \cdots, m\}$ of cardinality $n$.

Since $\mu = s_1 (f_1 + \sum_{j=1}^{m-n} \ga_j g_{1j}) = 0$,
spurious roots are the zeros of the system
\begin{equation}
\left\{ \
\begin{array}{l}
G^Ts = -s_0G^Tb, \\
s_i \left(\ga_0 f_i + \sum_{j=1}^{m-n} \ga_j g_{ij}\right) = 0, \ i=1,\cdots, m,\\
h_{\primaldual}(\ga_0,\ga,s_0,s) = 0 \\
\end{array}
\right.
\label{spurious}
\end{equation}

Let $(\ga_0,\ga,s_0,s)$ be a spurious root. Since spurious roots are
finite we assume that $\ga_0 = 1$ and $s_0=1$. We set
$K = K(s) = \{ k \in \{1, \cdots, m \}: s_k = 0 \}$. Then, the system
(\ref{spurious}) breaks into:
\begin{equation}
\left\{ \
\begin{array}{l}
G^T (s + b) = 0, \text{ with } s \in S_K\\
f_i + \sum_{j=1}^{m-n} \ga_j g_{ij} = 0, \ i \not \in K\\
h_{\primaldual}(\ga,s) = 0 \\
\end{array}
\right.
\label{spurious2}
\end{equation}

The first equation is equivalent to say that $s \in S_K$ and
$s+ b \in (\ker A^H)^{\perp}$. Again, this is the same as
saying that
$s \in S_K \cap \left( (\ker A^H)^{\perp} - b \right)$.

If the cardinality of
$K$ is strictly larger than $n$, then
Proposition~\ref{genericity} item~\ref{gen-primal} guarantees that
the intersection of $S_K$ with the affine space
$(\ker A^H)^{\perp} + b $ is empty. Therefore the cardinality of $K$
is at most $n$.

The second equation of (\ref{spurious2}) is satisfied if and only
if the preimage of $c$ by $A^T$ and $S_{\{1,\dots,m\} \setminus K}$
have a non-empty intersection. In that
case, $\gamma_i$ is associated to the $i$-th non-zero
coordinate of a point $y$ at the intersection.

If the cardinality of $K$ is strictly smaller than $n$, then
the intersection is empty because of Proposition~\ref{genericity}
item~\ref{gen-dual}.

Therefore $K$ has cardinality $n$. Reciprocally, let $K$ be
a cardinality $n$ subset of $\{1, \cdots, m \}$. Then
because of Proposition~\ref{genericity} item~\ref{gen-primal},
$S_K$ and $(\ker A^H)^{\perp} - b$ do intersect. Because
of item~\ref{gen-dual}, the spaces
$(A^T)^{-1}(c)$ and $S_{\{1,\dots,m\} \setminus K}$ also have a
non-empty intersection.

Therefore, equations $(\ref{spurious2}-1,2)$ admit a finite common solution
$(1,\ga,1,s)$.
Here, it is crucial to use the fact that we are in the primal-dual or
dual case, and hence $h_*$ is the determinant of $M''$.
It will follow
that $h_*(\ga_0,\ga,s_0,s) = 0$. This happens because
for $i \notin K$, $Q_{i}$ will vanish in (\ref{eq-mu-aug}).
Also for $i \notin K $, $M'(\ga,\mu)_{i, m+1} = 0$.

Since $\mu$ was replaced by
$s_1 (f_1 + \sum \ga_j g_{1j})$ in the last row of (\ref{eq-mu-aug3}),
the matrix  $M'(\ga,\mu)$ has (after reordering)
an $(m-n+1) \times (n+1)$ block of
zeros. This means that an $n+1$-dimensional subspace is mapped into
a $m+1-(m-n+1) = n$-dimensional
subspace. Hence $h_{*}(\ga,s) = \det M'(\ga,\mu) = 0$.

\end{proof}

Notice that spurious roots must be isolated. Otherwise, there would be
a component of non-isolated roots of $(\ref{E3})$ for $\mu \ne 0$, and
this would contradict Proposition~\ref{genericity} item~\ref{gen-nondeg}.
Hence, when $*$ is one of $\primaldual$ or $\dual$, we can subtract the number
of spurious roots from the bounds obtained
in Lemma~\ref{lem-main-th-9}, and Proposition~\ref{prop-mu} follows.
Notice that the bound for the dual case is not sharp.
}

\section{Concluding remarks.}

1. Beling and Verma \cite{bel} is a predecessor to our paper. They
prove a similar result to our Proposition~\ref{prop-mu} but only
for subspaces defined by the vanishing of one of the coordinates
and their estimate is not as strong.

\vskip 3mm 2. We have estimated the curvature by the number of
complex roots of a system of equations including possibly roots at
infinity. In fact only real and finite roots count. The number of
real roots is in general much less and can in some contexts be
compared with the square root of the number of complex roots, see
Shub and Smale \cite{Bez2}, Edelman and Kostlan \cite{ede},
McLennan \cite{mac}, Rojas \cite{roj} and Malajovich and Rojas
\cite{mal,mal2}. Thus the total curvature
at least on average may be very small indeed for large problems.
We find a better understanding of the total curvature of the
central path in worst and average case analysis an interesting
problem.

\vskip 3mm 3. There is a body of literature on the curvature of
the central path, relating the curvature to the complexity of
Newton type algorithms that approximate the central path and
produce approximations to the solutions: see Sonnevend, Stoer and
Zhao \cite{son1,son2}, Stoer and Zhao \cite{sto}, Zhao
\cite{zha1, zha2}. These papers use a different notion
of curvature, closer to $1 / \ga$ where $\ga$ is Smale's $\ga$,
see also Dedieu and Smale \cite{ded2}. The integral of these
quantities is infinite.

\vskip 3mm 4. The Riemannian geometry of the central path has been
studied by quite a few authors, see Karmarkar \cite{kar}, Bayer and Lagarias \cite{bay1a,bay1b,bay2}, Nesterov and Todd \cite{nes}.

\vskip 3mm 5. Vavasis and Ye \cite{vav} study regions where the
tangent vectors to the central path stay in definite cones.
Curvature estimates may be used as a refinement of this
information.

\vskip 3mm {6.
Malajovich and Meer~\cite{MM} showed that the problem of
computing (or even approximating up to
a fixed constant) the sharpest multi-homogeneous
B\'ezout bound for a system of polynomial equations is
actually $\mathrm{NP}$-hard.}

\vskip 3mm

{\bf Acknowledgements}

We would like to thank the two anonymous referees and the
communicating editor, Mike Todd, for their advice and for
suggesting that we state some of our results in terms of sign
invariant measures.

\providecommand{\bysame}{\leavevmode\hbox to3em{\hrulefill}\thinspace}

\end{document}